\UseRawInputEncoding
\documentclass[12pt,reqno]{amsart}
\usepackage{amsmath,amsfonts,amsthm,amsopn,amssymb}
\usepackage{cite,marginnote}
\usepackage{color,enumitem,graphicx}
\usepackage[colorlinks=true,urlcolor=blue,
citecolor=red,linkcolor=blue,linktocpage,pdfpagelabels,
bookmarksnumbered,bookmarksopen]{hyperref}
\usepackage[english]{babel}

\usepackage[left=2.9cm,right=2.9cm,top=2.8cm,bottom=2.8cm]{geometry}
\usepackage[hyperpageref]{backref}

\numberwithin{equation}{section}

\makeindex
\makeindex
\newtheorem{theorem}{Theorem}[section]
\newtheorem{definition}[theorem]{Definition}
\newtheorem{lemma}[theorem]{Lemma}
\newtheorem{corollary}[theorem]{Corollary}
\newtheorem{proposition}[theorem]{Proposition}
\newtheorem{remark}[theorem]{Remark}

\newcommand{\s}{\section}

\newcommand{\R}{\mathbb R}

\newcommand{\lab}{\label}
\newcommand{\bt}{\begin{theorem}}
\newcommand{\et}{\end{theorem}}
\newcommand{\bl}{\begin{lemma}}
\newcommand{\el}{\end{lemma}}
\newcommand{\bd}{\begin{definition}}
\newcommand{\ed}{\end{definition}}
\newcommand{\bc}{\begin{corollary}}
\newcommand{\ec}{\end{corollary}}
\newcommand{\bp}{\begin{proof}}
\newcommand{\ep}{\end{proof}}
\newcommand{\bx}{\begin{example}}
\newcommand{\ex}{\end{example}}
\newcommand{\bi}{\begin{exercise}}
\newcommand{\ei}{\end{exercise}}
\newcommand{\bo}{\begin{proposition}}
\newcommand{\eo}{\end{proposition}}
\newcommand{\br}{\begin{remark}}
\newcommand{\er}{\end{remark}}
\newcommand{\beq}{\begin{equation}}
\newcommand{\eeq}{\end{equation}}
\newcommand{\ba}{\begin{align}}
\newcommand{\ea}{\end{align}}
\newcommand{\bn}{\begin{enumerate}}
\newcommand{\en}{\end{enumerate}}
\newcommand{\bg}{\begin{align*}}
\newcommand{\bcs}{\begin{cases}}
\newcommand{\ecs}{\end{cases}}

\newcommand{\bean}{\begin{eqnarray*}}
\newcommand{\eean}{\end{eqnarray*}}

\def\R{\mathbb{R}}

\def\bd{\mathrm{bd}\,}

\title[Existence and Concentration]{Existence and Concentration Results
for the General Kirchhoff Type Equations}

\author[Y.~B.~Deng]{Yinbin Deng}
\author[W.~Shuai]{Wei Shuai}
\author[X.~X.~Zhong]{Xuexiu Zhong}

\address[Y.~B.~Deng]{\newline\indent ~School of Mathematics and Statistics \& Hubei Key Laboratory of Mathematical Sciences
\newline\indent
Central China Normal University
\newline\indent
 Wuhan 430079, P. R. China}
\email{\href{mailto:ybdeng@mail.ccnu.edu.cn}{ybdeng@ccnu.edu.cn}}

\address[W.~Shuai]{\newline\indent School of Mathematics and Statistics \& Hubei Key Laboratory of Mathematical Sciences
\newline\indent
Central China Normal University
\newline\indent
Wuhan 430079, P. R. China}
\email{\href{mailto:wshuai@mail.ccnu.edu.cn}{wshuai@ccnu.edu.cn}}

\address[X.~X.~Zhong]{\newline\indent South China Research Center for Applied Mathematics and Interdisciplinary Studies
\newline\indent
South China Normal University
\newline\indent
Guangzhou 510631, PR China}
\email{\href{mailto:zhongxuexiu1989@163.com}{zhongxuexiu1989@163.com}}

\thanks{The research is partially supported by the NSFC (No.11801581, 11931012,12071170), Guangdong Basic and Applied Basic Research Foundation (2021A1515010034), Guangzhou Basic and Applied Basic Research Foundation(202102020225) and Province Natural Science Fund of Guangdong (2018A030310082). }

\subjclass[2000]{}

\begin{document}

\begin{abstract}
We consider the following singularly perturbed  Kirchhoff type equations
$$-\varepsilon^2 M\left(\varepsilon^{2-N}\int_{\R^N}|\nabla u|^2 dx\right)\Delta u +V(x)u=|u|^{p-2}u~\hbox{in}~\R^N, u\in H^1(\R^N),N\geq 1,$$
where $M\in C([0,\infty))$ and $V\in C(\R^N)$ are given functions. Under very mild assumptions on $M$,  we prove the existence of single-peak or   multi-peak solution $u_\varepsilon$ for above problem,  concentrating around topologically stable critical points of $V$, by a direct corresponding argument. This gives an  affirmative answer to an open problem raised by Figueiredo et al. in 2014 [ARMA,213].

\end{abstract}

\maketitle

\s{Introduction}
\renewcommand{\theequation}{1.\arabic{equation}}

In the present paper, we study the existence and the concentration behavior  of positive solutions to the general Kirchhoff type equations
\beq\lab{eq:main-P}
-\varepsilon^2 M\left(\varepsilon^{2-N}\int_{\R^N}|\nabla u|^2 dx\right)\Delta u +V(x)u=f(u)~\hbox{in}~\R^N, u\in H^1(\R^N),u>0~\hbox{in}~\R^N,
\eeq
where $M\in C([0,\infty))$, $f\in C(\R)$ and $V\in C(\R^N)$ are given functions. If $M(t)\equiv 1$, it becomes the well-known nonlinear Schr\"odinger equation (replace $\varepsilon$ by $\delta$):
\beq\lab{eq:main-P2}
-\delta^2 \Delta w+V(x)w=f(w), w>0, w\in H^1(\R^N).
\eeq
In the past decades, a lot of work devote to the study of semiclassical solutions for \eqref{eq:main-P2}.

Recalling the  pioneering work \cite{Floer1986}, Floer and Weinstein  first studied the existence of single-peak solutions for $N=1,V\in C^2(\R^N)$ and $f(s)=s^3$. They construct a single-peak solution concentrating around any given nondegenerate critical point of $V$. And the higher dimension case with $f(s)=|s|^{p-2}s$, $p\in (2,2^*)$ is studied by Oh\cite{Oh1988}. In \cite{Floer1986,Oh1988}, their arguments are based on a Lyapunov-Schmidt reduction,  which requires a linearized nondegeneracy of a solution for a limiting problem. That is, if
$$-\Delta \phi +m\phi-f'(W)\phi=0~\hbox{in}~\R^N, \phi\in H^1(\R^N),$$
then $\phi=\sum_{i=1}^{N}a_i \frac{\partial W}{\partial x_i}$ for some $a_i\in \R$. Here $W$ is a ground state solution of the autonomous problem
 \beq\lab{eq:20220302-xe1}
    -\Delta w+mw=f(w), w\in H^1(\R^N).
 \eeq
Moreover, they also required that $f(t)/t$ is nondecreasing on $(0,\infty)$ and the uniqueness of ground states of \eqref{eq:20220302-xe1}. After then, many authors have applied Lyapunov-Schmidt reduction approach to  further refined results for more general $f$ and more general types of critical points of $V$ (see \cite{Ambrosetti1997,Ambrosetti2003,Pino2007,Kang2000,Li1997,Oh1990} and references therein).
It is known that the linearized nondegeneracy condition holds only for  a restricted class of $f$ for $N\geq 2$. We remark that  one needs at least the monotonicity of $\displaystyle \frac{mt-f'(t)t}{mt-f(t)}$ for $t>T$, where $T$ is the first positive zero of $mt-f(t)=0$ (see \cite{Cortazar1998}).  Even though there is such a restriction on the nonlinearity when one applies the reduction method, the Lyapunov-Schmidt reduction method is a very powerful tool when one constructs very subtle (highly unstable) solutions with continuum peaks as we can see in \cite{Pino2007}.
 Dancer developed a refined finite-dimensional reduction to construct peak solutions without the linearized nondegeneracy condition in \cite{Dancer2009}. However, he still
requires some type of nondegeneracy for the limiting problem.

We also remark that the variational approach is also proved very effective to study the problem \eqref{eq:main-P2}, which does not require the nondegeneracy condition for the limiting problem \eqref{eq:20220302-xe1}. This kind of approach was initiated by Rabinowitz\cite{Rabinowitz.1992} and has been developed further by several authors (see \cite{Byeon2003,Pino1996,Pino1997,Pino2002,Gui1996,Byeon2007,Byeon2008,Byeon2010,Byeon2013a,Jeanjean2004,Zhang2014,Zhang2015,Byeon2013,Zhang2017,Wang1993} and the references therein). In \cite{Rabinowitz.1992}, Rabinowitz proved the existence of positive solutions of \eqref{eq:main-P2} for small $\delta$ whenever $V\in C(\R^N,\R)$ and  $\displaystyle \liminf_{|x|\rightarrow \infty}V(x)>\inf_{x\in \R^N}V(x)=:V_0>0$.
Wang  proved that these solutions (obtained by Rabinowitz\cite{Rabinowitz.1992}) concentrate around the global minimum points of $V$ as $\delta\rightarrow 0$ in \cite{Wang1993}. In \cite{Pino1996},  del Pino and Felmer established a localized version of the result in \cite{Rabinowitz.1992,Wang1993}. Precisely, suppose the following
\begin{itemize}
\item[$(V_1)$] $V\in C(\R^N,\R)$ and $\displaystyle \inf_{x\in \R^N}V(x)=V_0>0$;
\item[$(V_2)$] there is a bounded domain $O$ such that
$$m:=\inf_{x\in O}V(x)<\inf_{x\in \partial O}V(x),$$
\end{itemize}
they obtained a single-peak solution concentrating around the minimum points of $V$ in $O$, provided that $f$ satisfies some conditions, such as the Ambrosetti-Rabinowitz condition and that the function $t\mapsto f(t)/t$ is nondecreasing.
After then, a lot of works are devoted to weak the assumptions on $f$ and construct peak solutions concentrating at more general critical points (such as local maximum points and special saddle points).

Byeon et al. \cite{Byeon2007,Byeon2008} developed a new variational method to explore what the essential features that guarantee the existence of localized ground states are. They studied the nonlinearities of Sobolev sub-critical case under the well know Berestycki-Lions conditions, which were first proposed in the pioneer work \cite{Berestycki1983} to guarantee the existence of ground states of \eqref{eq:20220302-xe1} in the subcritical case. So Byeon and Jeanjean \cite{Byeon2007} believed that Berestycki-Lions conditions are almost optimal for the subcritical case.
Byeon and Tanaka \cite{Byeon2013a} improved the result of \cite{Byeon2007,Byeon2008} by proving the existence of positive solutions of \eqref{eq:main-P2} also under the Berestycki-Lions conditions, which concentrate at more general critical points of $V$, such as saddle points or local maximum points.
We also remark that d'Avenia et al. \cite{Pietro2012} developed a min-max argument to establish the concentration phenomenon around the saddle points of the potential $V$.
In \cite{Zhang2014}, Zhang et al. generalized the result of \cite{Byeon2007} to the nonlinearities involving critical growth.
In \cite{Zhang2015}, Zhang and Zou also established the concentration phenomenon around the saddle points of the potential $V$ for nonlinearities involving critical growth, which generalized the result of subcritical case given by d'Avenia et al. \cite{Pietro2012}.

When $M(t)\not\equiv const$, He and Zou\cite{He2012} study \eqref{eq:main-P} when $N=3$ and $M(t)=a+bt$,  and show the existence of positive solutions of \eqref{eq:main-P} concentrating to global minima of $V$, under some suitable assumptions on $f$ which is of sub-critical case. Later, Wang et al. \cite{Wang2012} generalized the result to the nonlinearities involving Sobolev critical case. The authors in \cite{He2012,Wang2012} used the Nehari manifold method and thus the positive solution obtained is indeed has the least energy among all nontrivial solutions of \eqref{eq:main-P}. In \cite{He2016}, Yi He studied Problem \eqref{eq:main-P} with the nonlinearity involving Sobolev critical case when $N=3$ and $M(t)=a+bt$.
Some other related results we refer to \cite{Figueiredo2012,Figueiredo2013,He2014,Li2014,Li2014a} and references therein.

Inspired by the work of \cite{Byeon2007,Byeon2008,Byeon2008b},   Figueiredo et al. \cite{Figueiredo2014} studied problem \eqref{eq:main-P} for general $M$ by purely variational approach. Under suitable conditions on $M$ and Berestycki-Lions conditions on $f$, they constructed  a family of positive solutions $u_\varepsilon$ (for sufficiently small $\varepsilon$, and may not be least energy solution of \eqref{eq:main-P}) which concentrates at a local minimum of $V$ up to a subsequence. Precisely, for $M\in C([0,\infty))$, suppose $(M_1)$ when $N =1,2$ and $(M_1)-(M_5)$ when $N\geq 3$ below:
\begin{itemize}
\item[$(M_1)$]There exists $m_0>0$ such that $M(t)\geq m_0>0$ for any $t\geq 0$.
\item[$(M_2)$]Set $\widehat{M}(t):=\int_0^t M(s)ds$. Then there holds $\displaystyle \liminf_{t\rightarrow +\infty} \left\{\widehat{M}(t)-(1-2/N)M(t)t\right\}=+\infty$.
\item[$(M_3)$]$M(t)/t^{2/(N-2)}\rightarrow 0$ as $t\rightarrow +\infty$.
\item[$(M_4)$] The function $M$ is nondecreasing in $[0,\infty)$.
\item[$(M_5)$]) The function $t\mapsto M(t)/t^{2/(N-2)}$ is nonincreasing in $(0,\infty)$.
\end{itemize}
And assume that $f$ satisfies the following (f1)-(f4).
\begin{itemize}
\item[$(f_1)$] $f\in C(\R), f(s)=0$ for $s\leq 0$.
\item[$(f_2)$]$\displaystyle -\infty<\liminf_{s\rightarrow 0^+}\frac{f(s)}{s}\leq \limsup_{s\rightarrow 0^+}\frac{f(s)}{s}<V_0:=\inf_{x\in \R^N}V(x)$.
\item[$(f_3)$]When $N\geq 3$, $f(s)/s^{2^*-1}\rightarrow 0$ as $s\rightarrow +\infty$ and when $N=2$, $f(s)/e^{\alpha s^2}\rightarrow 0$ as $s\rightarrow +\infty$, for any $\alpha>0$, where $2^*:=2N/(N-2)$.
\item[$(f_4)$]There exists $T>0$ such that if $N\geq 2$, $(m/2)T^2<F(T)$ and if $N=1$, $\frac{1}{2}mt^2>F(t)$ for $t\in (0,T)$, $\frac{1}{2}mT^2=F(T)$ and $mT<f(T)$, where $F(t):\equiv\int_0^tf(s)ds$.
\end{itemize}
Then Figueiredo et al. established the following result.

\noindent
{\bf Theorem  A  (concentrating to local minimum points, see \cite[Theorem 1.1]{Figueiredo2014}).} \\
Assume $(V_1)-(V_2)$ and $(f_1)-(f_4)$. In addition, suppose $(M_1)$ when $N=1,2$ and $(M_1)-(M_5)$ when $N\geq 3$. Then there exists $\bar{\varepsilon}>0$ and a family $(u_\varepsilon)_{0<\varepsilon<\bar{\varepsilon}}$ of positive solutions of \eqref{eq:main-P} satisfying the following:
\begin{itemize}
\item[(i)] $dist(x_\varepsilon, \mathcal{M})\rightarrow 0$ as  $\varepsilon\rightarrow 0$, where $(x_\varepsilon)$ be a maximum point of $u_\varepsilon$ and $\mathcal{M}\equiv \left\{x\in O: V(x)=m\right\}$.
\item[(ii)] After taking a subsequence $(\varepsilon_n), u_{\varepsilon_n}(\varepsilon_n x+x_{\varepsilon_n})\rightarrow U$ strongly in $H^1(\R^N)$, where $U$ is a positive least energy solution of
    \beq\lab{eq:20220303-e2}
    -M\left(\int_{\R^N}|\nabla u|^2 dx\right)\Delta u+mu=f(u)~\hbox{in}~\R^N, u\in H^1(\R^N).
    \eeq
\item[(iii)] There exists $C_1,C_2>0$ such that
\beq\lab{eq:20220303-e3}
u_\varepsilon (x)\leq C_1 \exp{\left(-C_2\frac{|x-x_\varepsilon|}{\varepsilon}\right)}~\hbox{for all $x\in \R^N$ and $0<\varepsilon<\bar{\varepsilon}$}.
\eeq
\end{itemize}\hfill$\Box$

\br\lab{remark:20220324-r1}
\    \
\begin{itemize}
\item[(i)]
If we study problem \eqref{eq:main-P} by purely variational method, the difficulty  is the presence of nonlocal term $M\left(\int_{\R^N}|\nabla u|^2 dx\right)$, which makes \eqref{eq:main-P} more  delicate than \eqref{eq:main-P2}.
\item[(ii)]Also if we study problem \eqref{eq:main-P} by Lyapunov-Schmidt reduction directly, one need the nondegenerate condition that the kernel of the  linearized operator $\mathcal{L}:L^2(\R^N)\mapsto L^2(\R^N)$  difinded by
\begin{align}\lab{eq:20220324-e2}
\mathcal{L}\varphi=&-M\left(\int_{\R^N}|\nabla U|^2 dx\right)\Delta \varphi-2M'\left(\int_{\R^N}|\nabla U|^2 dx\right)\left(\int_{\R^N}\nabla U\cdot \nabla \varphi dx\right)\Delta U\nonumber\\
& \quad+m\varphi -f'(U)\varphi
\end{align}
is spanned by the functions $\frac{\partial U}{\partial x_i}, i=1,\cdots,N$, provided $M,f\in C^1$.
 We remark that if $N=3, M(t)=a+bt$ with $a>0,b>0$ and $f(u)=|u|^{p-1}u, 1<p<5$, Li et al.\cite{Li2020} proved that there exists a unique positive radial solution $U\in H^1(\R^3)$ satisfying
\beq\lab{eq:20220324-e1}
-\left(a+b\int_{\R^3}|\nabla U|^2 dx\right)\Delta U+U=U^p, U>0~\hbox{in}~\R^3.
\eeq
Moreover $U$ is nondegenerate in $H^1(\R^3)$ (see \cite[Theorem 1.2]{Li2020}). Then the authors in \cite{Li2020} could apply the Lyapunov-Schmidt reduction to construct a family of solutions concentrating at a local minimum (see \cite[Theorem 1.3]{Li2020}). Furthermore, they obtain the local uniqueness of single peak solutions provided some further  assumptions on $V(x)$ (see \cite[Theorem 1.4]{Li2020}). However, for the general $M(t)$, the non-degenerate condition is very hard to check even for the nonlinearities are of polynomials. Hence, it is also very hard to study problem \eqref{eq:main-P} by a direct Lyapunov-Schmidt reduction.
\end{itemize}
\er

So Figueiredo et al. raised the following {\bf open problem} (see \cite[Remark 1.2-(iii)]{Figueiredo2014}):
\begin{itemize}
  \item [ \ ]{\it ``It seems interesting to consider whether one can find a family of solutions
of \eqref{eq:main-P} which has multi-peaks or which concentrates around other type of
critical points of $V$ (local maxima, saddle points and so on). These types of results to \eqref{eq:main-P2} have been obtained." }
\end{itemize}


There is little progress on this open problem.
As far as we know, in \cite{Chen2019}, Chen and Ding gave a first affirmative answer to this open problem for $N\geq 3$, by constructing positive solutions concentrating around the local maximum points of $V$, basing on the same assumptions on $M$ and $f$ as in \cite{Figueiredo2014} and in addition that $M(t)+(1-N/2)M'(t)t\neq 0$.

In present paper, we shall give another affirmative answer to this open problem, by constructing positive single-peak solutions concentrating around local minimum, local maxima or saddle points. Our assumptions on $M$ are mild. Furthermore, we also obtain some result about the concentrating positive multi-peak solutions, which seem never be obtained in the related literatures.

The paper is organized as follows. In the next section, we give the corresponding relationship  between \eqref{eq:main-P} and \eqref{eq:main-P2} for single-peak solutions and multi-peak solutions respectively. Taking $f(s)=|s|^{p-2}s$ with $2<p<2^*$ as applications, we obtain the existence and multiplicity of single-peak solutions  (see Theorem \ref{th:main-th1}) and the existence of multi-peak solutions (see Theorem \ref{th:20220326-xth1}) under different assumptions on $V(x)$ and mild assumption on $M(t)$.
In the section \ref{sec:proof-5}, we shall prove these corresponding theorems. In the section \ref{sec:proof-1-4}, we will give some sufficient conditions that guarantee the corresponding theorems to be applied. Theorem \ref{th:main-th1} and Theorem \ref{th:20220326-xth1} will be proved in Section \ref{sec:proof-7}.

\s{Statement of main results}\lab{sec:statement}
\renewcommand{\theequation}{2.\arabic{equation}}

We define
$$H_\varepsilon=\left\{u(x)\in H^1(\R^N), \int_{\R^N}\big(\varepsilon^2 |\nabla u(x)|^2 +V(x)u^2(x)\big)dx<\infty\right\},$$
and for any $u(x)\in H_\varepsilon$, denote its norm by
\beq\lab{eq:def-norm}
\|u\|_\varepsilon:= (u(x), u(x))_{\varepsilon}^{\frac{1}{2}}=\left(\int_{\R^N}\big(\varepsilon^2 |\nabla u(x)|^2 +V(x)u^2(x)\big)dx\right)^{\frac{1}{2}}.
\eeq
For $u\in L^p(\R^N), 1\leq p<\infty$, we denote the $L^p$-norm by $\|\cdot\|_p$ for simplicity, i.e.,
$$\|u\|_p:=\left(\int_{\R^N}|u|^p dx\right)^{\frac{1}{p}}.$$

\subsection{Some correspondent results}
As the reasons stated  in Remark \ref{remark:20220324-r1}, we will not deal with \eqref{eq:main-P} directly. We firstly establish the following corresponding results between \eqref{eq:main-P} and \eqref{eq:main-P2} for single-peak solutions and multi-peak solutions respectively.

\bt\lab{th:main-t5}{\bf (correspondence for single-peak solution)}
Suppose, under some conditions on $V$ and $f$, that there exists $\bar{\delta}>0$ and a family $(\omega_\delta)_{0<\delta<\bar{\delta}}$ of positive solutions of \eqref{eq:main-P2} satisfying the following:
\begin{itemize}
\item[$(a$-$i)$] $dist(x_\delta, \mathcal{M})\rightarrow 0$ as  $\delta\rightarrow 0$, where $\mathcal{M}$ is an isolated set of topologically stable critical points of $V$ and $(x_\delta)$ is a maximum point of $\omega_\delta$;
\item[$(a$-$ii)$] after taking a subsequence $(\delta_n), \omega_{\delta_n}(\delta_n x+x_{\delta_n})\rightarrow W$ strongly in $H^1(\R^N)$, where $W$ is a positive least energy solution of
    \beq\lab{eq:20220303-xe1}
    -\Delta w+mw=f(w)~\hbox{in}~\R^N, w\in H^1(\R^N),
    \eeq
 $m:=V(x_0)$ with some $x_0\in \mathcal{M}$ and $x_{\delta_n}\rightarrow x_0$ as $n\rightarrow +\infty$;
\item[$(a$-$iii)$] there exists $C_1,C_2>0$ such that
\beq\lab{eq:20220303-xe2}
\omega_\delta (x)\leq C_1 \exp{\left(-C_2\frac{|x-x_\delta|}{\delta}\right)}~\hbox{for all $x\in \R^N$ and $0<\delta<\bar{\delta}$}.
\eeq
\end{itemize}
In addition,  suppose, under some conditions on $M$, that
\begin{itemize}
\item[$(b$-$i)$] there exists $\bar{\varepsilon}>0$ such that for any $\varepsilon\in (0,\bar{\varepsilon})$, there exists $\delta_\varepsilon\in (0, \bar{\delta})$ and a positive solution $\omega_{\delta_\varepsilon}$ of \eqref{eq:main-P2} such that
    \beq\lab{eq:20220303-xe3}
    \varepsilon^2 M\left(\varepsilon^{2-N}\|\nabla \omega_{\delta_\varepsilon}\|_{L^2}^{2}\right)=\delta_\varepsilon^2;
    \eeq
\item[(b-ii)] it holds that
\beq\lab{eq:20220303-xe4}
0<\inf_{\varepsilon\in (0,\bar{\varepsilon})}\frac{\delta_\varepsilon}{\varepsilon}\leq \sup_{\varepsilon\in (0,\bar{\varepsilon})}\frac{\delta_\varepsilon}{\varepsilon}<+\infty.
\eeq
\end{itemize}
For $\varepsilon\in (0,\bar{\varepsilon})$, define $u_\varepsilon(x):\equiv \omega_{\delta_\varepsilon}(x)$ and denote a maximum point $x_{\delta_\varepsilon}$ of $\omega_{\delta_\varepsilon}$ by $x_\varepsilon$ for simplicity. Then $(u_\varepsilon)_{0<\varepsilon<\bar{\varepsilon}}$ is a family of positive solutions to \eqref{eq:main-P} and $x_\varepsilon$ is also a maximum point of $u_\varepsilon$. Furthermore,
\begin{itemize}
\item[$(c$-$i)$] $dist(x_\varepsilon, \mathcal{M})\rightarrow 0$ as  $\varepsilon\rightarrow 0$;
\item[$(c$-$ii)$] after taking a subsequence $(\varepsilon_n), u_{\varepsilon_n}(\varepsilon_n x+x_{\varepsilon_n})\rightarrow U$ strongly in $H^1(\R^N)$, where $U$ is a positive least energy solution of
    \beq\lab{eq:20220303-xe5}
    -M(\|\nabla w\|_{L^2}^{2})\Delta w+mw=f(w)~\hbox{in}~\R^N, w\in H^1(\R^N),
    \eeq
 $m:=V(x_0)$ with some $x_0\in \mathcal{M}$ and $x_{\varepsilon_n}\rightarrow x_0$ as $n\rightarrow +\infty$;
\item[$(c$-$iii)$] there exists $C_3,C_4>0$ such that
\beq\lab{eq:20220303-xe6}
u_\varepsilon (x)\leq C_3 \exp{\left(-C_4\frac{|x-x_\varepsilon|}{\varepsilon}\right)}~\hbox{for all $x\in \R^N$ and $0<\varepsilon<\bar{\varepsilon}$}.
\eeq
\end{itemize}
\et

\bt\lab{th:main-bt5}{\bf (correspondence for multi-peak solution)}
Suppose, under some conditions on $V$ and $f$, that there exists $\bar{\delta}>0$ and a family $(\omega_\delta)_{0<\delta<\bar{\delta}}$ of multi-peak positive solutions of \eqref{eq:main-P2} such that $\omega_\delta$ is the form of
\beq\lab{eq:20220325-xe1}
\omega_{\delta}(x)=\sum_{j=1}^{k} W_{{\bf P}_j}\left(\frac{x-y_{\delta}^{(j)}}{\delta}\right)+\psi_\delta(x)
\eeq
with $y_{\delta}^{(j)},\psi_\delta(x)$ satisfying
\beq\lab{eq:20220325-xe2}
y_{\delta}^{(j)}\rightarrow {\bf P}_j~\hbox{and}~ \|\psi_\delta\|_\delta=o(\delta^{\frac{N}{2}})
\eeq
as  $\delta\rightarrow 0$ for $j=1,\cdots,k$,
where $W_{{\bf P}_j}$ is a radial positive ground state solution of
\beq\lab{eq:20220325-xe3}
-\Delta w+V({\bf P}_j) w=f(w)~\hbox{in}~\R^N, w\in H^1(\R^N).
\eeq
and $P_j,  \ j=1,2, \cdots, k$, are critical points of $V(x)$.
In addition, suppose, under some conditions on $M$, that $M$ satisfies the assumptions $(b$-$i)$ and (b-ii) in Theorem \ref{th:main-t5}.

For $\varepsilon\in (0,\bar{\varepsilon})$, define $u_\varepsilon(x):= \omega_{\delta_\varepsilon}(x), \phi_\varepsilon(x):=\psi_{\delta_\varepsilon}(x)$ and $ x_{\varepsilon}^{(j)}=y_{\delta_\varepsilon}^{(j)}$. Then $(u_\varepsilon)_{0<\varepsilon<\bar{\varepsilon}}$ is a family of multi-peak positive solutions to \eqref{eq:main-P}. In particular, $u_\varepsilon$ is the form of
\beq\lab{eq:20220326-e1}
u_\varepsilon(x)=\sum_{j=1}^{k}U_{{\bf P}_j}\left(\frac{x-x_{\varepsilon}^{(j)}}{\varepsilon}\right)+\phi_\varepsilon(x),
\eeq
with $x_{\varepsilon}^{(j)}, \phi_\varepsilon(x)$ satisfying
\beq\lab{eq:20220326-e2}
x_{\varepsilon}^{(j)}\rightarrow {\bf P}_j~\hbox{and}~\|\phi_\varepsilon\|_{\varepsilon}=o(\varepsilon^{\frac{N}{2}}).
\eeq
as $\varepsilon\rightarrow 0$ for $j=1,\cdots,k$.
Here $(U_{{\bf P}_1},\cdots,U_{{\bf P}_k})$ is a positive solution to the following system
\beq\lab{eq:20220325-xe7}
-M\left(\sum_{j=1}^{k}\|\nabla u_j\|_2^2\right)\Delta u_j+V({\bf P}_j)u_j=f(u_j), j=1,\cdots,k.
\eeq
\et

\br\lab{remark:20220325-r1}
 \    \
\begin{itemize}
\item[(i)] Indeed, up to a subsequence, we may assume that $\frac{\delta_\varepsilon}{\varepsilon}\rightarrow C_*$ as $\varepsilon\rightarrow 0$. Then $U_{{\bf P}_j}(x):=W_{{\bf P}_j}\left(\frac{x}{C_*}\right)$. In particular, $C_*$ solves $G(t)=0$ with $G(t)$ defined by
    \beq\lab{eq:20220326-e3}
    G(t):=M\left(t^{N-2}\big(\sum_{j=1}^{k}\|\nabla W_{{\bf P}_j}\|_2^2\big)\right)-t^2,\ \ \  t>0.
    \eeq
\item[(ii)]Under some conditions on $M$, if $G(t)=0$ has a unique positive root, one can see that $(U_{{\bf P}_1},\cdots,U_{{\bf P}_k})$ is indeed the unique positive radial solution to the system \eqref{eq:20220325-xe7}. For example, if $M(t)=a+bt, a,b>0,N=3$, one can easily verify that
    $$G(t)=a+b\big(\sum_{j=1}^{k}\|\nabla W_{{\bf P}_j}\|_2^2\big)t-t^2.$$
    It is trivial that $G(t)=0$ has a unique positive root. In such a case, our result coincides \cite[Proposition 3]{Luo2019}, where the authors deal with \eqref{eq:main-P} by Lyapunov-Schmidt reduction method directly, thanks to the nondegeneracy of positive solutions to the limit Kirchhoff problem proved in \cite{Li2020}.
\item[(iii)] If $G(t)=0$ has multiple positive roots, it will be interesting to consider the multiplicity of multi-peak solutions concentrating to the same given $k$ critical points of $V(x)$.
\end{itemize}
\er

\subsection{Single-peak solutions: Existence and multiplicity}

As applications, we firstly concerned with the existence and multiplicity of single-peak solutions of \eqref{eq:main-P}.
We adopt some definitions by Grossi\cite{Grossi2002}.
\begin{definition}\lab{def:d1}
We say that a function $h:\R^N\mapsto \R$ is homogenous of degree $\alpha\in \R^+$ with respect to ${\bf P}\in \R^N$ if
\beq\lab{eq:20220324-xe1}
h(tx+{\bf P})=t^\alpha h(x+{\bf P})~\hbox{for any $t\in \R^+$ and $x\in \R^N$}.
\eeq
\end{definition}

\begin{definition}\lab{def:d2}
{\bf(Definition of admissible potential)} Let us assume that $V\in C^1(\R^N)$ satisfies
\beq\lab{eq:20220324-xe2}
|\nabla V(x)|\leq C e^{\gamma |x|} ~\hbox{at infinity}
\eeq
and
\beq\lab{eq:20220324-xe3}
0<V_0\leq V(x)\leq V_1
\eeq
for some $\gamma>0$. We say that $V$ is an admissible potential at ${\bf P}\in \R^N$ if there exist continuous functions $h_i:\R^N\mapsto \R, \ \  R_i: B_{{\bf P},r}\equiv\left\{x\in \R^N: |x-{\bf P}|<r\right\}\mapsto \R$ and real numbers $\alpha_i\geq 1, i=1,\cdots,N$, such that
\begin{itemize}
\item[(i)] $\displaystyle \frac{\partial V}{\partial x_i}(x)=h_i(x)+R_i(x)$ in $B_{{\bf P},r}$;
\item[(ii)] $R_i(x)\leq C |x-{\bf P}|^{\beta_i}$ in $B_{{\bf P},r}$ with $\beta_i>\alpha_i$ for any $i=1,\cdots,N$;
\item[(iii)] $h_i(x)=0$ if and only if $x={\bf P}$;
\item[(iv)] $h_i$ is homogeneous of degree $\alpha_i$ respect to  {\bf P}.
\end{itemize}
\end{definition}

\begin{definition}\lab{def:d3}
Let $G\in C(\R^N,\R^N)$ be a vector field. We say that $y$ is a stable zero for $G$ if
\begin{itemize}
\item[(i)] $G(y)=0$,
\item[(ii)] $y$  isolated,
\item[(iii)] if $G_n$ is a sequence of vector fields such that $\|G_n-G\|_{C(B_{y,\rho})}\rightarrow 0$ as $n \rightarrow \infty$ for some $\rho>0$, then there exists $y_n$ such that $G_n(y_n)=0$ and $y_n\rightarrow y$ as $n \rightarrow \infty$.
\end{itemize}
\end{definition}

Set the following vector field
\beq\lab{eq:20220324-xe4}
\mathcal{L}_{\bf P}(y)=\left(\int_{\R^N}h_i(x+y+{\bf P})W_{\bf P}^{2}(x)\right)_{i=1,\cdots,N},
\eeq
where $W_{\bf P}\in H^1(\R^N)$ is the unique positive radial solution to
\beq\lab{eq:20220515-e1}
-\Delta w+V({\bf P})w=w^{p-1}~\hbox{in}~\R^N.
\eeq
Define
\beq\lab{eq:20220324-xe5}
Z=\{y\in \R^N~\hbox{such that $y$ is a stable zero of $\mathcal{L_{\bf P}}$}\},
\eeq
then we get the following theorem which is concerned with the existence and multiplicity of single-peak solution concentrating at general critical point ${\bf P}$.
\bt\lab{th:main-th1}
Let $V$ be an admissible potential at $P$. Suppose that $\# Z<\infty$ and
\beq\lab{eq:20220324-xe6}
det~Jac~\mathcal{L}_{\bf P}(y)\neq 0
\eeq
for any $y\in Z$.
In addition, suppose $(M_1)$ when $N=1,2$ and suppose $(M_1)$ and $(M_3)$ when $N\geq 3$.
Then there exists $\bar{\varepsilon}>0$ such that for any $0<\varepsilon<\bar{\varepsilon}$, the equation
 \beq\lab{eq:20220324-xe7}
 \begin{cases}
-\varepsilon^2 M\left(\varepsilon^{2-N}\int_{\R^N}|\nabla u|^2 dx\right)\Delta u +V(x)u=|u|^{p-2}u~\hbox{in}~\R^N,\\
 u\in H^1(\R^N),u>0~\hbox{in}~\R^N, 2<p<2^*:=\frac{2N}{(N-2)_+}
 \end{cases}
\eeq
possesses at least $\#Z$ different single-peak solutions concentrating at $x={\bf P}$.
Precisely, we obtain families of positive solutions $(u_{\varepsilon}^{(k)})_{0<\varepsilon<\bar{\varepsilon}, k=1,\cdots,\#Z}$ with the maximum point at
$({\bf P}+\varepsilon y_{\varepsilon}^{(k)})_{0<\varepsilon<\bar{\varepsilon}, k=1,\cdots,\#Z}$
such that, for any $k\in \{1,\cdots,\#Z\}$, $y_{\varepsilon}^{(k)}$ is bounded with
$y_{\varepsilon}^{(k)}\rightarrow y_k$ as $\varepsilon\rightarrow 0$. Here $\{y_k, k=1,\cdots,\#Z\}$ are $\#Z$ distinct stable zeros of $\mathcal{L}_{\bf P}(y)$. In particular,
$u_{\varepsilon}^{(k)}(\varepsilon x+{\bf P}+\varepsilon y_{\varepsilon}^{(k)})\rightarrow U$ strongly in $H^1(\R^N)$,  where $U$ is a positive radial ground state solution of
\beq\lab{eq:20220324-xe8}
-M\left(\|\nabla u\|_2^2\right)\Delta u+V({\bf P})u=|u|^{p-2}u~\hbox{in}~\R^N.
\eeq
\et

When ${\bf P}$ is a nondegenerate critical point of $V$,
a direct conclusion of Theorem \ref{th:main-th1} can be stated as follows.

\bc\lab{cro:20220324-c1}{\bf (single-peak solution concentrating at nondegenerate critical point)}
Let ${\bf P}$ be a nondegenerate critical point of $V$ and $0<V_0\leq V(x)\leq V_1$. Suppose $(M_1)$ when $N=1,2$ and suppose $(M_1)$ and $(M_3)$ when $N\geq 3$. Then problem (\ref {eq:20220324-xe7}) possesses a family of positive solutions $(u_\varepsilon)_{0<\varepsilon<\bar{\varepsilon}}$ with the  maximum point at  $({\bf P}+\varepsilon y_\varepsilon)_{0<\varepsilon<\bar{\varepsilon}}$, such that $y_\varepsilon\rightarrow 0$ and $u_\varepsilon (\varepsilon x+{\bf P}+\varepsilon y_\varepsilon)\rightarrow U$ strongly in $H^1(\R^N)$ as $\varepsilon\rightarrow 0^+$, where $U$ is a positive radial ground state solution of \eqref{eq:20220324-xe8}.
\ec

\subsection{Multi-peak solutions}

The existence of multi-peak solutions for nonlinear Schr\"odinger equations \eqref{eq:main-P2} (or its variants) concentrating at the critical points of $V(x)$ also has been studied deeply, we refer to \cite{Pino1996,DelPino1998,Cao1996,Li1998,Noussair2000}.
For the case of a critical nonlinearity, the results on the existence of multi-peak solutions can be seen in \cite{Bahri1995,Rey1990}. For the case of existence with the concerntration phenomena, we refer to \cite{Ambrosetti2001,Ambrosetti2007,Cao2009,Cingolani2000,Pino1997,Deng2014,Gui1996,Kang2000} and the references therein. For the uniqueness of multi-bump solution, we refer to \cite{Cao2003,Cao2015}.
There is also some results focused on the multi-peak solutions for \eqref{eq:main-P} when $M(t)=a+bt, a>0,b>0$, see \cite{Luo2019,Wang2020}.

To state our result on the multi-peak solutions, we consider a class of $V(x)$ as follows:
\begin{itemize}
\item[$(V'_1)$] $V(x)\in C^1(\R^N)$ and $\inf_{x\in \R^N}=:V_0>0$.
\item[$(V'_2)$] $V(x)$ satisfies the following expansions:
\beq\lab{eq:20220325-ze1}
\begin{cases}
V(x)=V({\bf P}_j)+\sum_{i=1}^{N} b_{j,i} |x_i- {\bf P}_{j,i}|^\alpha +O(|x-{\bf P}_j|^{\alpha+1}), \quad & x\in B_{\eta} ({\bf P}_j),\\
\frac{\partial V(x)}{\partial x_i}=\alpha b_{j,i} |x_i- {\bf P}_{j,i}|^{\alpha-2}(x_i- {\bf P}_{j,i})+O(|x-{\bf P}_j|^{\alpha}), \quad & x\in B_{\eta} ({\bf P}_j),
\end{cases}
\eeq
where $\eta>0$ is a small constant, $\alpha>1, \  x=(x_1,\cdots,x_N),$  $ {\bf P}_j=({\bf P}_{j,1},\cdots,{\bf P}_{j,N}),$  $ b_{j,i}\in \R$ with $b_{j,i}\neq 0$ for each $i=1,\cdots,N, j=1,\cdots,k$.
\end{itemize}
Here comes to our main result about multi-peak solutions.
\bt\lab{th:20220326-xth1}
Assume that $V(x)$ satisfies $(V'_1)$ and $(V'_2)$. Suppose $(M_1)$ when $N=1,2$ and suppose $(M_1)$ and $(M_3)$ when $N\geq 3$. Then there exists $\bar{\varepsilon}>0$ and a family of positive solutions $(u_\varepsilon)_{0<\varepsilon<\bar{\varepsilon}}$ of \eqref{eq:20220324-xe7} concentrating at a set of $k$ different points $\{{\bf P}_1,\cdots,{\bf P}_k\}\subset \R^N$. Precisely,  $u_\varepsilon$ is the form of
\beq
u_\varepsilon(x)=\sum_{j=1}^{k}W_{{\bf P}_j}\left(\frac{x-x_{\varepsilon}^{(j)}}{\varepsilon}\frac{1}{C_*}\right)+\phi_\varepsilon(x),
\eeq
with $x_{\varepsilon}^{(j)},\phi_\varepsilon(x)$ satisfying
\beq\lab{eq:20220325-ze2}
|x_{\varepsilon}^{(k)}-{\bf P}_j|=o(\varepsilon)~\hbox{and}~\|\phi_\varepsilon\|_\varepsilon=O(\varepsilon^{\frac{N}{2}+\alpha})
\eeq
as $\varepsilon\rightarrow 0$ for $j=1,\cdots k$.

Here $W_{{\bf P}_j}$ is a radial positive ground state solution of \eqref {eq:20220325-xe3} with $f(w)=w^{p-1}$
and $C_*$ is the smallest positive number determined by
\beq\lab{eq:20220325-bze2}
M\left(C_{*}^{N-2}\big(\sum_{j=1}^{k}\|\nabla W_{{\bf P}_j}\|_2^2\big)\right)=C_{*}^{2}.
\eeq
\et

\s{Proofs of Theorem \ref{th:main-t5} and Theorem \ref{th:main-bt5}}\lab{sec:proof-5}
\renewcommand{\theequation}{3.\arabic{equation}}
\noindent{\bf Proof of Theorem \ref{th:main-t5}.}
By the assumption $(b$-$i)$  and the definition of $u_\varepsilon$, a direct computation shows that
\begin{align*}
&-\varepsilon^2 M\left(\varepsilon^{2-N}\|\nabla u_\varepsilon\|_2^2\right)\Delta u_\varepsilon +V(x)u_\varepsilon-f(u_\varepsilon)\\
=&-\varepsilon^2 M\left(\varepsilon^{2-N}\|\nabla \omega_{\delta_\varepsilon}\|_2^2\right)\Delta \omega_{\delta_\varepsilon} +V(x)\omega_{\delta_\varepsilon}-f(\omega_{\delta_\varepsilon})\\
=&-\delta_\varepsilon^2 \Delta \omega_{\delta_\varepsilon} +V(x)\omega_{\delta_\varepsilon}-f(\omega_{\delta_\varepsilon})\\
=&0.
\end{align*}
Hence, $u_\varepsilon$ is a solution to \eqref{eq:main-P} and it is clear that $x_\varepsilon$ is a maximum point of $u_\varepsilon$. Thus the conclusion of $(c$-$i)$ holds by $(a$-$i)$.

By the assumption $(b$-$ii)$, we can find some $K_1>0$ such that
\beq\lab{eq:20220303-xe7}
\frac{\delta_\varepsilon}{\varepsilon}\leq K_1, \forall \varepsilon\in (0,\bar{\varepsilon}).
\eeq
Noting that $\delta_\varepsilon<\bar{\delta}, \varepsilon\in (0,\bar{\varepsilon})$, by the assumption $(a$-$iii)$ and the definition of $u_\varepsilon$, we have that
\begin{align*}
u_\varepsilon(x)=\omega_{\delta_\varepsilon}(x)\leq&C_1 \exp{\left(-C_2\frac{|x-x_{\delta_\varepsilon}|}{\delta_\varepsilon}\right)}\\
=&C_1 \exp{\left(-C_2\frac{|x-x_{\varepsilon}|}{\varepsilon}\cdot \frac{\varepsilon}{\delta_\varepsilon}\right)}\\
\leq&C_1 \exp{\left(-\frac{C_2}{K_1}\frac{|x-x_{\varepsilon}|}{\varepsilon}\right)}~\hbox{for all $x\in \R^N$ and $0<\varepsilon<\bar{\varepsilon}$}.
\end{align*}
Hence, the conclusion $(c$-$iii)$ holds.

After taking a subsequence $(\varepsilon_n)$, by the assumption (b-ii), we may assume that
\beq\lab{eq:20220303-xe8}
\frac{\delta_{\varepsilon_n}}{\varepsilon_n}\rightarrow C_*\in (0,K_1].
\eeq
Noting that $\delta_{\varepsilon_n}\rightarrow 0$, by the assumption of $(a$-$ii)$, we have that
$\varphi_n(x):\equiv\omega_{\delta_{\varepsilon_n}}(\delta_{\varepsilon_n} x+x_{\delta_{\varepsilon_n}})\rightarrow W$ strongly in $H^1(\R^N)$, where $W$ is a positive least energy solution of \eqref{eq:20220303-xe1}.
Hence,
\beq\lab{eq:20220303-xe9}
\|\nabla u_{\varepsilon_n}\|_2^2=\|\nabla \omega_{\delta_{\varepsilon_n}}\|_2^2
=\delta_{\varepsilon_n}^{N-2}\|\nabla \varphi_n\|_2^2=\delta_{\varepsilon_n}^{N-2}(\|\nabla W\|_2^2+o(1)).
\eeq
We note that \eqref{eq:20220303-xe3} in the assumption $(b$-$i)$ implies that
\beq\lab{eq:20220303-xe10}
M\left(\varepsilon_{n}^{2-N}\|\nabla u_{\varepsilon_n}\|_2^2\right)=\left(\frac{\delta_{\varepsilon_n}}{\varepsilon_n}\right)^2.
\eeq
Let $n\rightarrow +\infty$, it follows from $M\in C([0,+\infty))$  that
\beq\lab{eq:20220303-xe11}
M\left({C_*}^{N-2}\|\nabla W\|_2^2\right)={C_*}^2.
\eeq
Put $U(x):\equiv W(\frac{1}{C_*}x)$, a direct computation shows that $U$ is a positive least energy solution to
\beq\lab{eq:20220303-xe12}
-C_*^2 \Delta u+mu=f(u)~\hbox{in}~\R^N, u\in H^1(\R^N).
\eeq
Noting that $\|\nabla U\|_2^2={C_*}^{N-2}\|\nabla W\|_2^2$, so by \eqref{eq:20220303-xe11} and \eqref{eq:20220303-xe12}, we see that $U$ is a positive least energy solution to
\eqref{eq:20220303-xe5}.
In particular,
\begin{align*}
&u_{\varepsilon_n}(\varepsilon_n x+x_{\varepsilon_n})=\omega_{\delta_{\varepsilon_n}}(\varepsilon_n x+x_{\varepsilon_n})\\
=&\omega_{\delta_{\varepsilon_n}}(\delta_{\varepsilon_n}\cdot \frac{\varepsilon_n}{\delta_{\varepsilon_n}} x+x_{\delta_{\varepsilon_n}})\\
=&\omega_{\delta_{\varepsilon_n}} \left(\delta_{\varepsilon_n} \left(\frac{1}{C^*}+o(1)\right)x+x_{\delta_{\varepsilon_n}}\right)\\
\rightarrow &W\left(\frac{1}{C^*}x\right)=U(x)~\hbox{in}~H^1(\R^N).
\end{align*}
Hence, the conclusion of $(c$-$ii)$ holds. We complete the proof of Theorem \ref{th:main-t5}.\hfill$\Box$

\vskip 0.2in
\noindent{\bf Proof of Theorem \ref{th:main-bt5}.}
Similar to the proof of Theorem \ref{th:main-t5}. We only note that in this case,
$$\|\nabla u_{\varepsilon_n}\|_2^2=\|\nabla \omega_{\delta_{\varepsilon_n}}\|_2^2
=\delta_{\varepsilon_n}^{N-2}\left(\sum_{j=1}^{k}\|\nabla W_{{\bf P}_j}\|_2^2+o(1)\right)$$
and thus $C_*$ satisfies
$$M\left({C_*}^{N-2}\big(\sum_{j=1}^{k}\|\nabla W_{{\bf P}_j}\|_2^2\big)\right)={C_*}^2,$$
i.e., $G(C_*)=0$.
We also note that $o(\delta_{\varepsilon}^{\frac{N}{2}})=o(\varepsilon^{\frac{N}{2}})$ due to $(b$-$ii)$.
\qed

\s{Some sufficient conditions to guarantee $(b$-$i)$ and $(b$-$ii)$}\lab{sec:proof-1-4}
\renewcommand{\theequation}{4.\arabic{equation}}

\bl\lab{lemma:20220303-xl1}
Under suitable assumptions on $V$ and $f$, we assume that $\omega_\delta$ depends continuously on $\delta$ and there exists some $A>0$ such that $\frac{\|\nabla \omega_\delta\|_2^2}{\delta^{N-2}}\rightarrow A$ as $\delta\rightarrow 0$.
Suppose $(M_1)$ when $N=1,2$ and suppose $(M_1)$ and $(M_3)$ when $N\geq 3$. Then the conditions $(b$-$i)$ and $(b$-$ii)$ in Theorem \ref{th:main-t5} hold.
\el
\bp
For $\delta>0$ small, and $\varepsilon>0$ small, we define
\beq\lab{eq:20220302-e1}
g_\varepsilon(\delta):=\varepsilon^2 M\left(\varepsilon^{2-N}\|\nabla \omega_\delta\|_2^2\right)-\delta^2.
\eeq
Since $\omega_\delta$ depends continuously on $\delta$  and from the assumption $(M_1)$, we see that
$g_\varepsilon (\delta)$ is continuous with respect to $\delta\in (0,\bar{\delta})$.
Let us consider the equation $g_\varepsilon(\delta)=0$.
By $(M_1)$, there exists some $K_0>0$ large enough such that
\beq\lab{eq:20220315-e0}
\begin{cases}
\frac{1}{K^2}M\left(\frac{A}{K}\right)<1, \ \ \forall K\geq K_0,\quad&\hbox{if}~N=1,\\
\frac{1}{K^2}M\left(A\right)<1, \ \ \forall K\geq K_0,\quad& \hbox{if}~N=2.
\end{cases}
\eeq
Suppose further $(M_3)$ if $N\geq 3$, for $K$ large enough, $(M_3)$ implies that
$$\frac{1}{K^2}M\left(K^{N-2}A\right)=\frac{1}{K^2} o_K\left(\big(K^{N-2}A\big)^{\frac{2}{N-2}}\right)=o_K(1).$$
So there also exists $K_0>0$ large enough such that
\beq\lab{eq:20220315-bue0}
\frac{1}{K^2}M\left(K^{N-2}A\right)<1, \ \ \forall K\geq K_0.
\eeq

Now, let $K=K_0$ be fixed.
We claim that there exists  $\delta_1<\bar{\delta}$ small enough such that
\beq\lab{eq:20220315-e1}
g_{\frac{\delta}{K}}(\delta)<0, \ \  \forall \delta\in (0,\delta_1).
\eeq
If not, there exists $\delta_n\rightarrow 0$ such that
\beq\lab{eq:20220315-e2}
g_{\frac{\delta_n}{K}}(\delta_n)\geq 0.
\eeq
Under the assumption, we have that
\beq\lab{eq:20220315-e4}
\|\nabla \omega_{\delta_n}\|_2^2=\delta_{n}^{N-2}(A+o(1)).
\eeq
Hence, by \eqref{eq:20220315-bue0}, for $n$ large enough,

\begin{align} \lab{eq:20220315-e5}
&\left(\frac{\delta_n}{K}\right)^2 M\left(\left(\frac{\delta_n}{K}\right)^{2-N} \delta_{n}^{N-2}(A+o(1))\right)-\delta_n^2  \nonumber  \\
&=\left(\frac{1}{K^2}M\left(K^{N-2}\big(A+o(1)\big)\right)-1\right)\delta_n^2<0,
\end{align}
a contraction. Hence, \eqref{eq:20220315-e1} holds. In another word, there exists $\varepsilon_1:=\frac{\delta_1}{K}>0$ small enough such that
\beq\lab{eq:20220315-e6}
g_{\varepsilon}(K\varepsilon)<0, \forall \varepsilon\in (0,\varepsilon_1).
\eeq
On the other hand, under the assumption $(M_1)$, we have that
\beq\lab{eq:20220303-xe13}
g_\varepsilon(\delta)\geq m_0 \varepsilon^2 -\delta^2>0~\hbox{for all} ~\delta\in (0,\sqrt{m_0}\varepsilon).
\eeq
Define
\beq\lab{eq:20220302-e2}
\delta_\varepsilon:=\sup\left\{s: g_\varepsilon(\delta)>0, \ \ \forall 0<\delta<s\right\}, \varepsilon\in (0,\varepsilon_1).
\eeq
By \eqref{eq:20220303-xe13} and \eqref{eq:20220315-e6}, one can see that $\delta_\varepsilon$ is well defined and thus $g_\varepsilon(\delta_\varepsilon)=0$. In particular,
\beq\lab{eq:20220315-bule1}
\sqrt{m_0}\varepsilon\leq\delta_\varepsilon\leq K\varepsilon, \ \ \forall \varepsilon\in (0,\varepsilon_1).
\eeq
Hence, $(b$-$i)$ and $(b$-$ii)$ hold.
\ep

\br\lab{remark:20220303-wr1}
\   \
\begin{itemize}
\item[(i)]
We remark that the assumption $(M_3)$ plays an crucial role to guarantee the existence of $\delta_\varepsilon$. This condition has requirements for dimension, and usually the high-dimensional case is not applicable.  Indeed, it is a sufficient but not necessary condition. For example, we take $M(t)=a+bt$ with $a,b>0$. For the case of $N\geq 4$, $(M_3)$ fails for any $a,b>0$. However, if $a$ and $b$ are small suitable such that
\beq\lab{eq:20220311-e1}
\inf\left\{M(A t^{N-2})-t^2, t>0\right\}=\inf\left\{a+bA t^{N-2}-t^2: t>0\right\}<0,
\eeq
then one can prove the existence of $\delta_\varepsilon$,  arguing by contradiction like the case $N\geq 3$ in the proof of Lemma \ref{lemma:20220303-xl1}. Then we can also establish the similar result for $M(t)=a+bt$ with $N\geq 4$, which is very difficult to obtain by a direct variational method.
\item[(ii)]Under some suitable assumptions on the general nonlinearity $f$, if $\inf_{x\in \R^N}V(x)=:V_0>0$ large enough, then the uniqueness and nondegeneracy of positive solution to \eqref{eq:20220325-xe3} hold, see \cite[Theorem 1.3]{JeanZhangZhong2021} and \cite{Cortazar1998}. In particular, $\|\nabla W_{{\bf P}_j}\|_2\rightarrow +\infty$ as $V({\bf P}_j)\rightarrow +\infty$.   Hence, $A\rightarrow \infty$ as $V({\bf P}_j)\rightarrow \infty$ and thus \eqref{eq:20220311-e1} is not expected provided $\inf_{x\in \R^N}V(x)=:V_0>0$ large enough. Indeed, if these kinds of concentrating results are valid independent of $V_0$, the condition $(M_3)$ is almost necessary for $N\geq 3$, see the following Proposition \ref{prop:20220311-p1}.
\end{itemize}
\er

\bo\lab{prop:20220311-p1}
Let $N\geq 3$. $M\in C([0,\infty))$ satisfies $(M_1)$ and
\beq\lab{eq:20220311-e2}
\liminf_{t\rightarrow +\infty}M(t)/t^{2/(N-2)}>0.
\eeq
Assume that $f\in C(\R,\R)$ satisfies the following conditions:
\begin{itemize}
\item[(i)]$\displaystyle \lim_{s\rightarrow 0}\frac{f(s)}{s}=0$;
\item[(ii)]there exists some $\ell\in (2,2^*)$ such that $\lim_{|s|\rightarrow \infty}\frac{|f(s)|}{|s|^{\ell-1}}<\infty$.
\end{itemize}
Then there exists some $V_0>0$ such that if $\inf_{x\in \R^N}V(x)\geq V_0$, the equation \eqref{eq:main-P} has no nontrivial solution in $H^1(\R^N)$ for all $\varepsilon>0$.
\eo
\bp
Noting that $u(x)$ solves \eqref{eq:main-P} if and only if $u(\varepsilon x)$ solves
\beq\lab{eq:20220311-e3}
-M(\|\nabla \phi \|_2^2)\Delta \phi +V(\varepsilon x)\phi=f(\phi)~\hbox{in}~\R^N, \phi\in H^1(\R^N),
\eeq
we only need to prove the conclusion holds for $\varepsilon=1$.
By $(M_1)$ and \eqref{eq:20220311-e2}, it is easy to see that
\beq\lab{eq:20220311-e4}
\inf_{t\geq 0}\frac{M(t)}{t^{\frac{2}{N-2}}}=:\sigma>0.
\eeq
Thus,
\beq\lab{eq:20220311-e5}
M(\|\nabla u\|_2^2)\|\nabla u\|_2^2\geq \sigma \|\nabla u\|_{2}^{2^*}, \forall u\in H^1(\R^N).
\eeq
For any $u$ solves
\beq\lab{eq:20220311-e6}
-M(\|\nabla u\|_2^2)\Delta u+V(x)u=f(u)~\hbox{in}~\R^N, u\in H^1(\R^N),
\eeq
we have that
\beq\lab{eq:20220311-e7}
M(\|\nabla u\|_2^2)\|\nabla u\|_2^2+\int_{\R^N}V(x)u^2dx =\int_{\R^N}f(u)u dx.
\eeq
We remark that under the assumptions on $f$, one can see that for any $\eta>0$, there exists some $C_\eta>0$ such that
\beq\lab{eq:20220311-e8}
|f(s)s|\leq \eta|s|^2+C_\eta |s|^{\ell}, \forall s\in \R.
\eeq
Hence, by \eqref{eq:20220311-e5}- \eqref{eq:20220311-e8}, we obtain that
\beq\lab{eq:20220311-e9}
\sigma \|\nabla u\|_{2}^{2^*}+(V_0-\eta) \|u\|_2^2\leq C_\eta \|u\|_{\ell}^{\ell}.
\eeq
Recalling the Gagliardo-Nierenberg inequality, there exists some $C_\ell>0$ such that
\beq\lab{eq:20220310-e6}
\|u\|_\ell^\ell \leq C_\ell \|\nabla u\|_{2}^{\frac{N(\ell-2)}{2}} \|u\|_{2}^{\ell-\frac{N(\ell-2)}{2}}.
\eeq
Thus, by \eqref{eq:20220311-e9} and \eqref{eq:20220310-e6}, there exists some $C^*=C_{\eta,\ell}>0$ such that
\beq\lab{eq:20220311-e10}
\sigma \|\nabla u\|_{2}^{2^*}+(V_0-\eta) \|u\|_2^2\leq C^*\|\nabla u\|_{2}^{\frac{N(\ell-2)}{2}} \|u\|_{2}^{\ell-\frac{N(\ell-2)}{2}}.
\eeq
By Young inequality, take $p=\frac{4}{(N-2)(\ell-2)}$ and $p'=\frac{4}{2\ell-N(\ell-2)}$, we have that
\beq\lab{eq:20220311-e11}
\|\nabla u\|_{2}^{\frac{N(\ell-2)}{2}} \|u\|_{2}^{\ell-\frac{N(\ell-2)}{2}}
\leq \kappa \|\nabla u\|_{2}^{2^*} + (p-1)p^{-\frac{p}{p-1}}\kappa^{-\frac{1}{p-1}} \|u\|_2^2, \forall \kappa>0.
\eeq
By taking $\kappa=\sigma$, it follows \eqref{eq:20220311-e10} and \eqref{eq:20220311-e11} that
\beq\lab{eq:20220311-e12}
(V_0-\eta) \|u\|_2^2\leq  C^*(p-1)p^{-\frac{p}{p-1}}\sigma^{-\frac{1}{p-1}}\|u\|_2^2, \hbox{where}~p=\frac{4}{(N-2)(\ell-2)}.
\eeq
Hence, if $V_0>\eta+C^*(p-1)p^{-\frac{p}{p-1}}\sigma^{-\frac{1}{p-1}}$ with $p=\frac{4}{(N-2)(\ell-2)}$, then $u\equiv 0$.
\ep

\s    {Proof of Theorem \ref{th:main-th1}  and Theorem \ref{th:20220326-xth1} }\lab{sec:proof-7}
\renewcommand{\theequation}{5.\arabic{equation}}
\noindent
{\bf Proof of Theorem \ref{th:main-th1}.}
We firstly recall some known results.
Since $V(x)$ is an admissible potential, by \cite[Theorem 4.3]{Grossi2002}, there exists $\bar{\delta}>0$ such that for any stable zero $y_k$ of $\mathcal{L}_{\bf P}(y)$,  one can construct single-peak solutions $(\omega_{\delta}^{(k)})_{0<\delta<\bar{\delta}}$ to equation \eqref{eq:main-P2} such that :
\begin{itemize}
\item[(i)]let $x_{\delta}^{(k)}$ be the maximum point of $\omega_\delta$, then $x_{\delta}^{(k)}$ can be written as
\beq
x_{\delta}^{(k)}={\bf P}+\delta y_{\delta}^{(k)}~\hbox{with}~y_{\delta}^{(k)}\rightarrow y_k~\hbox{as}~\delta\rightarrow 0.
\eeq
\item[(ii)] $\omega_\delta(\delta x+x_{\delta}^{(k)})\rightarrow W_{\bf P}$  in $H^1(\R^N)$.
\item[(iii)] For $j\neq k$, the single-peak solutions generated by $y_j$ and $y_k$ are different.
\end{itemize}

Furthermore, if $det~Jac~\mathcal{L}_{\bf P}(y_k)\neq 0$, then $y_{\delta}^{(k)}$ is unique determined in a small neighborhood of $y_k$. Hence, \eqref{eq:main-P2} possesses exactly $\#Z$ single-peak solutions concentrating at ${\bf P}$, see \cite[Theorem 1.1]{Grossi2002}.

So for any fixed $k$ and $\delta>0$ small, the uniqueness of $y_{\delta}^{(k)}$ implies that $(\omega_{\delta}^{(k)})_{0<\delta<\bar{\delta}}$ depends on $\delta$ continuously.
Indeed, we can take $\eta>0$ small enough such that
$B_\eta(y_j)\cap B_\eta(y_k)=\emptyset, \forall j\neq k$. Take $\bar{\delta}$ small enough such that $y_{\delta}^{(k)}\in B_\eta(y_k)$ for $\forall \delta\in (0,\bar{\delta})$ and $k=1,\cdots,\#Z$.
Let $k$ be fixed, for any $\delta^*\in (0,\bar{\delta})$ and any sequence $\delta_n\rightarrow \delta^*$. Up to a subsequence, we may assume that $y_{\delta_n}^{(k)}\rightarrow y_{*}^{(k)}\in \overline{B_\eta(y_k)}$. It is trivial that
\beq\lab{eq:20220326-xe2}
W_{\bf P}\left(\frac{x-({\bf P}+\delta_n y_{\delta_n}^{(k)})}{\delta_n}\right)\rightarrow W_{\bf P}\left(\frac{x-({\bf P}+\delta^* y_{*}^{(k)})}{\delta^*}\right)~\hbox{strongly in}~H^1(\R^N).
\eeq
Define


\beq\lab{eq:20220326-xe1}
E_{\delta,{\bf P},y}=\left\{w(x)\in H^1(\R^N):\ \ \left . \begin {cases} &\left(w(x),W_{\bf P}\left(\frac{x-y}{\delta}\right)\right)_\delta=0, \\
&\left(w(x),\frac{W_{\bf P}\left(\frac{x-y}{\delta}\right)}{\partial x_i}\right)_\delta=0, i=1,\cdots,N
\end {cases} \right \}
\right\}.
\eeq

The standard Lyapunov-Schmidt reduction process implies that $\omega_{\delta}^{(k)}$ is of form
\beq\lab{eq:20220326-xe2}
\omega_{\delta}^{(k)}(x)=W_{\bf P}\left(\frac{x-x_{\delta}^{(k)}}{\delta}\right)+\psi_{\delta}^{(k)}(x),
\eeq
with $\psi_{\delta}^{(k)}\in E_{\delta,{\bf P},x_{\delta}^{(k)}}^{\perp}$.
By the uniformly decay of $\{\omega_{\delta_n}^{(k)}\}$, we may assume that
\beq\lab{eq:20220326-xe3}
\omega_{\delta_n}^{(k)}\rightarrow \omega^*~\hbox{strongly in}~H^1(\R^N).
\eeq
So that $\omega^*$ is a solution to
\beq\lab{eq:20220326-xe4}
-{\delta^*}^2\Delta w+V({\bf P})w =|w|^{p-2}w~\hbox{in}~\R^N, 0<w\in H^1(\R^N)
\eeq
 Furthermore, it is of form
\beq\lab{eq:20220326-xe5}
\omega^*=W_{\bf P}\left(\frac{x-({\bf P}+\delta^* y_{*}^{(k)})}{\delta^*}\right)+\psi^*,
\eeq
with $\psi_{\delta_n}^{(k)}\rightarrow \psi^*$ strongly in $H^1(\R^N)$.
Hence, one can see that $\psi^*\in E_{\delta,{\bf P},{\bf P}+\delta^* y_{*}^{(k)}}^{\perp}$.
Hence, by the uniqueness of $y_{\delta}^{(k)}$, we have that $y_{*}^{(k)}=y_{\delta^*}^{(k)}$ and thus $\omega^*=\omega_{\delta^*}^{(k)}$. Hence, $(\omega_{\delta}^{(k)})$ is continuous with  respect to $\delta\in (0,\bar{\delta})$.

On the other hand, by $\omega_\delta(\delta x+x_{\delta}^{(k)})\rightarrow W_{\bf P}$  in $H^1(\R^N)$, it is trivial that
$$\frac{\|\nabla \omega_\delta\|_2^2}{\delta^{N-2}}\rightarrow A:=\|\nabla W_{\bf P}\|_2^2~\hbox{as}~\delta\rightarrow 0^+.$$
So combining with Lemma \ref{lemma:20220303-xl1}, we see that Theorem \ref{th:main-t5} is applicable here and we finish the proof of Theorem \ref{th:main-th1}.\hfill$\Box$

\vskip 0.2in
\noindent
{\bf Proof of Corollary \ref{cro:20220324-c1}.} It is a direct conclusion of Theorem \ref{th:main-th1}. We only note that when ${\bf P}$ is a nondegenerate critical point of $V(x)$, one has that $Z=\{0\}$ and thus $\#Z=1$.\hfill$\Box$

\renewcommand{\theequation}{6.\arabic{equation}}

\noindent
{\bf Proof of Theorem \ref{th:20220326-xth1}.}
Firstly, we remark that under the assumptions, for the problem
\beq
-\delta^2 \Delta w+V(x)w=|w|^{p-2}w~\hbox{in}~\R^N, w\in H^1(\R^N),
\eeq
 the existence of multi-peak solutions concentrating at critical points of $V(x)$ is standard by Lyapunov-Schmidt reduction method, see for example \cite{Cao-book}. Secondly, under the assumptions of Theorem \ref{th:20220326-xth1},  Cao et al. also proved the local uniqueness result, see \cite[Theorem 1.1]{Cao2015}. Precisely, for $\delta$ small enough, $\omega_\delta(x)$ is of form
 \beq\lab{eq:20220326-we1}
 \omega_{\delta}(x)=\sum_{j=1}^{k}W_{{\bf P}_j}\left(\frac{x-y_{\delta}^{(j)}}{\delta}\right)+\psi_\delta(x)
 \eeq
 with $y_{\delta}^{(j)},\psi_\delta(x)$ satisfying, for $j=1,\cdots,k$, as $\delta\rightarrow 0$,
 \beq\lab{eq:20220326-we2}
 |y_{\delta}^{(j)}-{\bf P}_j|=o(\delta)~\hbox{and}~\|\psi_\delta\|_\delta=O(\delta^{\frac{N}{2}+\alpha})=o(\delta^{\frac{N}{2}}).
 \eeq
By $\delta^2\|\nabla \psi_\delta\|_2^2\leq \|\psi_\delta\|_\delta^2$, we see that
$\|\nabla \psi_\delta\|_2^2=o(\delta^{N-2})$. Hence, it follows that
$$\frac{\|\nabla \omega_\delta\|_2^2}{\delta^{N-2}}\rightarrow A:=\sum_{j=1}^{k}\|\nabla W_{{\bf P}_j}\|_2^2~\hbox{as}~\delta\rightarrow 0^+.$$
We  also remark that the local uniqueness can imply that $(\omega_{\delta})_{0<\delta<\bar{\delta}}$ depends continuously on $\delta$, see also the proof of Theorem  \ref{th:main-th1}. So the conclusions of Theorem \ref{th:20220326-xth1} follow by Theorem \ref{th:main-bt5}. We only note that by \eqref{eq:20220326-we2} and that
 $$0<\liminf_{\varepsilon\rightarrow 0}\frac{\delta_\varepsilon}{\varepsilon}\leq \limsup_{\varepsilon\rightarrow 0}\frac{\delta_\varepsilon}{\varepsilon}<+\infty,$$
 we obtain
 $$\|\phi_\varepsilon\|_\varepsilon=\|\psi_{\delta_\varepsilon}\|_\varepsilon
 =O(\|\psi_{\delta_\varepsilon}\|_{\delta_\varepsilon})=O(\delta_{\varepsilon}^{\frac{N}{2}+\alpha})
 =O(\varepsilon^{\frac{N}{2}+\alpha}).$$\hfill$\Box$

 \noindent
{\bf Acknowledgements}\\
The authors thanks Prof. Peng Luo for the valuable discussions when preparing the paper.


\begin{thebibliography}{10}

\bibitem{Ambrosetti1997}
A.~Ambrosetti, M.~Badiale, and S.~Cingolani.
\newblock Semiclassical states of nonlinear {S}chr\"{o}dinger equations.
\newblock {\em Arch. Rational Mech. Anal.}, 140(3):285--300, 1997.

\bibitem{Ambrosetti2001}
A.~Ambrosetti, A.~Malchiodi, and S.~Secchi.
\newblock Multiplicity results for some nonlinear {S}chr\"{o}dinger equations
  with potentials.
\newblock {\em Arch. Ration. Mech. Anal.}, 159(3):253--271, 2001.

\bibitem{Ambrosetti2007}
Antonio Ambrosetti and Andrea Malchiodi.
\newblock Concentration phenomena for nonlinear {S}chr\"{o}dinger equations:
  recent results and new perspectives.
\newblock In {\em Perspectives in nonlinear partial differential equations},
  volume 446 of {\em Contemp. Math.}, pages 19--30. Amer. Math. Soc.,
  Providence, RI, 2007.

\bibitem{Ambrosetti2003}
Antonio Ambrosetti, Andrea Malchiodi, and Wei-Ming Ni.
\newblock Singularly perturbed elliptic equations with symmetry: existence of
  solutions concentrating on spheres. {I}.
\newblock {\em Comm. Math. Phys.}, 235(3):427--466, 2003.

\bibitem{Bahri1995}
Abbas Bahri, Yanyan Li, and Olivier Rey.
\newblock On a variational problem with lack of compactness: the topological
  effect of the critical points at infinity.
\newblock {\em Calc. Var. Partial Differential Equations}, 3(1):67--93, 1995.

\bibitem{Berestycki1983}
H.~Berestycki and P.-L. Lions.
\newblock Nonlinear scalar field equations. {I}. {E}xistence of a ground state.
\newblock {\em Arch. Rational Mech. Anal.}, 82(4):313--345, 1983.

\bibitem{Byeon2010}
Jaeyoung Byeon.
\newblock Singularly perturbed nonlinear {D}irichlet problems with a general
  nonlinearity.
\newblock {\em Trans. Amer. Math. Soc.}, 362(4):1981--2001, 2010.

\bibitem{Byeon2007}
Jaeyoung Byeon and Louis Jeanjean.
\newblock Standing waves for nonlinear {S}chr\"{o}dinger equations with a
  general nonlinearity.
\newblock {\em Arch. Ration. Mech. Anal.}, 185(2):185--200, 2007.

\bibitem{Byeon2008b}
Jaeyoung Byeon and Louis Jeanjean.
\newblock Erratum: ``{S}tanding waves for nonlinear {S}chr\"{o}dinger equations
  with a general nonlinearity'' [{A}rch. {R}ation. {M}ech. {A}nal. {\bf 185}
  (2007), no. 2, 185--200; mr2317788].
\newblock {\em Arch. Ration. Mech. Anal.}, 190(3):549--551, 2008.

\bibitem{Byeon2008}
Jaeyoung Byeon, Louis Jeanjean, and Kazunaga Tanaka.
\newblock Standing waves for nonlinear {S}chr\"{o}dinger equations with a
  general nonlinearity: one and two dimensional cases.
\newblock {\em Comm. Partial Differential Equations}, 33(4-6):1113--1136, 2008.

\bibitem{Byeon2013a}
Jaeyoung Byeon and Kazunaga Tanaka.
\newblock Semi-classical standing waves for nonlinear {S}chr\"{o}dinger
  equations at structurally stable critical points of the potential.
\newblock {\em J. Eur. Math. Soc. (JEMS)}, 15(5):1859--1899, 2013.

\bibitem{Byeon2003}
Jaeyoung Byeon and Zhi-Qiang Wang.
\newblock Standing waves with a critical frequency for nonlinear
  {S}chr\"{o}dinger equations. {II}.
\newblock {\em Calc. Var. Partial Differential Equations}, 18(2):207--219,
  2003.

\bibitem{Byeon2013}
Jaeyoung Byeon, Jianjun Zhang, and Wenming Zou.
\newblock Singularly perturbed nonlinear {D}irichlet problems involving
  critical growth.
\newblock {\em Calc. Var. Partial Differential Equations}, 47(1-2):65--85,
  2013.

\bibitem{Cao1996}
Daomin Cao, Norman~E. Dancer, Ezzat~S. Noussair, and Shunsen Yan.
\newblock On the existence and profile of multi-peaked solutions to singularly
  perturbed semilinear {D}irichlet problems.
\newblock {\em Discrete Contin. Dynam. Systems}, 2(2):221--236, 1996.

\bibitem{Cao2003}
Daomin Cao and Hans-Peter Heinz.
\newblock Uniqueness of positive multi-lump bound states of nonlinear
  {S}chr\"{o}dinger equations.
\newblock {\em Math. Z.}, 243(3):599--642, 2003.

\bibitem{Cao2015}
Daomin Cao, Shuanglong Li, and Peng Luo.
\newblock Uniqueness of positive bound states with multi-bump for nonlinear
  {S}chr\"{o}dinger equations.
\newblock {\em Calc. Var. Partial Differential Equations}, 54(4):4037--4063,
  2015.

\bibitem{Cao2009}
Daomin Cao and Shuangjie Peng.
\newblock Semi-classical bound states for {S}chr\"{o}dinger equations with
  potentials vanishing or unbounded at infinity.
\newblock {\em Comm. Partial Differential Equations}, 34(10-12):1566--1591,
  2009.

\bibitem{Cao-book}
Daomin Cao, Shuangjie Peng, and Shusen Yan.
\newblock {\em Singularly perturbed methods for nonlinear elliptic problems}.
\newblock Cambridge Univ Press.

\bibitem{Chen2019}
Yu~Chen and Yanheng Ding.
\newblock Multiplicity and concentration for {K}irchhoff type equations around
  topologically critical points in potential.
\newblock {\em Topol. Methods Nonlinear Anal.}, 53(1):183--223, 2019.

\bibitem{Cingolani2000}
Silvia Cingolani and Monica Lazzo.
\newblock Multiple positive solutions to nonlinear {S}chr\"{o}dinger equations
  with competing potential functions.
\newblock {\em J. Differential Equations}, 160(1):118--138, 2000.

\bibitem{Cortazar1998}
Carmen Cort\'{a}zar, Manuel Elgueta, and Patricio Felmer.
\newblock Uniqueness of positive solutions of {$\Delta u+f(u)=0$} in {${\bf
  R}^N,\ N\ge3$}.
\newblock {\em Arch. Rational Mech. Anal.}, 142(2):127--141, 1998.

\bibitem{Dancer2009}
E.~N. Dancer.
\newblock Peak solutions without non-degeneracy conditions.
\newblock {\em J. Differential Equations}, 246(8):3077--3088, 2009.

\bibitem{Pietro2012}
Pietro d'Avenia, Alessio Pomponio, and David Ruiz.
\newblock Semiclassical states for the nonlinear {S}chr\"{o}dinger equation on
  saddle points of the potential via variational methods.
\newblock {\em J. Funct. Anal.}, 262(10):4600--4633, 2012.

\bibitem{Pino2002}
Manuel del Pino and Patricio Felmer.
\newblock Semi-classical states of nonlinear {S}chr\"{o}dinger equations: a
  variational reduction method.
\newblock {\em Math. Ann.}, 324(1):1--32, 2002.

\bibitem{Pino1996}
Manuel del Pino and Patricio~L. Felmer.
\newblock Local mountain passes for semilinear elliptic problems in unbounded
  domains.
\newblock {\em Calc. Var. Partial Differential Equations}, 4(2):121--137, 1996.

\bibitem{Pino1997}
Manuel del Pino and Patricio~L. Felmer.
\newblock Semi-classical states for nonlinear {S}chr\"{o}dinger equations.
\newblock {\em J. Funct. Anal.}, 149(1):245--265, 1997.

\bibitem{DelPino1998}
Manuel Del~Pino and Patricio~L. Felmer.
\newblock Multi-peak bound states for nonlinear {S}chr\"{o}dinger equations.
\newblock {\em Ann. Inst. H. Poincar\'{e} Anal. Non Lin\'{e}aire},
  15(2):127--149, 1998.

\bibitem{Pino2007}
Manuel del Pino, Michal Kowalczyk, and Jun-Cheng Wei.
\newblock Concentration on curves for nonlinear {S}chr\"{o}dinger equations.
\newblock {\em Comm. Pure Appl. Math.}, 60(1):113--146, 2007.

\bibitem{Deng2014}
Yinbin Deng, Shuangjie Peng, and Huirong Pi.
\newblock Bound states with clustered peaks for nonlinear {S}chr\"{o}dinger
  equations with compactly supported potentials.
\newblock {\em Adv. Nonlinear Stud.}, 14(2):463--481, 2014.

\bibitem{Figueiredo2013}
Giovany~M. Figueiredo.
\newblock Existence of a positive solution for a {K}irchhoff problem type with
  critical growth via truncation argument.
\newblock {\em J. Math. Anal. Appl.}, 401(2):706--713, 2013.

\bibitem{Figueiredo2014}
Giovany~M. Figueiredo, Norihisa Ikoma, and Jo\~{a}o~R. Santos~J\'{u}nior.
\newblock Existence and concentration result for the {K}irchhoff type equations
  with general nonlinearities.
\newblock {\em Arch. Ration. Mech. Anal.}, 213(3):931--979, 2014.

\bibitem{Figueiredo2012}
Giovany~M. Figueiredo and Jo\~{a}o~R. Santos, Junior.
\newblock Multiplicity of solutions for a {K}irchhoff equation with subcritical
  or critical growth.
\newblock {\em Differential Integral Equations}, 25(9-10):853--868, 2012.

\bibitem{Floer1986}
Andreas Floer and Alan Weinstein.
\newblock Nonspreading wave packets for the cubic {S}chr\"{o}dinger equation
  with a bounded potential.
\newblock {\em J. Funct. Anal.}, 69(3):397--408, 1986.

\bibitem{Grossi2002}
Massimo Grossi.
\newblock On the number of single-peak solutions of the nonlinear
  {S}chr\"{o}dinger equation.
\newblock {\em Ann. Inst. H. Poincar\'{e} Anal. Non Lin\'{e}aire},
  19(3):261--280, 2002.

\bibitem{Gui1996}
Changfeng Gui.
\newblock Existence of multi-bump solutions for nonlinear {S}chr\"{o}dinger
  equations via variational method.
\newblock {\em Comm. Partial Differential Equations}, 21(5-6):787--820, 1996.

\bibitem{He2012}
Xiaoming He and Wenming Zou.
\newblock Existence and concentration behavior of positive solutions for a
  {K}irchhoff equation in {$\Bbb R^3$}.
\newblock {\em J. Differential Equations}, 252(2):1813--1834, 2012.

\bibitem{He2016}
Yi~He.
\newblock Concentrating bounded states for a class of singularly perturbed
  {K}irchhoff type equations with a general nonlinearity.
\newblock {\em J. Differential Equations}, 261(11):6178--6220, 2016.

\bibitem{He2014}
Yi~He, Gongbao Li, and Shuangjie Peng.
\newblock Concentrating bound states for {K}irchhoff type problems in {$\Bbb
  R^3$} involving critical {S}obolev exponents.
\newblock {\em Adv. Nonlinear Stud.}, 14(2):483--510, 2014.

\bibitem{Jeanjean2004}
Louis Jeanjean and Kazunaga Tanaka.
\newblock Singularly perturbed elliptic problems with superlinear or
  asymptotically linear nonlinearities.
\newblock {\em Calc. Var. Partial Differential Equations}, 21(3):287--318,
  2004.

\bibitem{JeanZhangZhong2021}
Louis Jeanjean, Jianjun Zhang, and Xuexiu Zhong.
\newblock A global branch approach to normalized solutions for the
  {S}chr\"odinger equation.
\newblock 2021.

\bibitem{Kang2000}
Xiaosong Kang and Juncheng Wei.
\newblock On interacting bumps of semi-classical states of nonlinear
  {S}chr\"{o}dinger equations.
\newblock {\em Adv. Differential Equations}, 5(7-9):899--928, 2000.

\bibitem{Li2020}
Gongbao Li, Peng Luo, Shuangjie Peng, Chunhua Wang, and Chang-Lin Xiang.
\newblock A singularly perturbed {K}irchhoff problem revisited.
\newblock {\em J. Differential Equations}, 268(2):541--589, 2020.

\bibitem{Li2014}
Gongbao Li and Hongyu Ye.
\newblock Existence of positive ground state solutions for the nonlinear
  {K}irchhoff type equations in {$\Bbb R^3$}.
\newblock {\em J. Differential Equations}, 257(2):566--600, 2014.

\bibitem{Li2014a}
Gongbao Li and Hongyu Ye.
\newblock Existence of positive solutions for nonlinear {K}irchhoff type
  problems in {$\Bbb R^3$} with critical {S}obolev exponent.
\newblock {\em Math. Methods Appl. Sci.}, 37(16):2570--2584, 2014.

\bibitem{Li1997}
YanYan Li.
\newblock On a singularly perturbed elliptic equation.
\newblock {\em Adv. Differential Equations}, 2(6):955--980, 1997.

\bibitem{Li1998}
Yanyan Li and Louis Nirenberg.
\newblock The {D}irichlet problem for singularly perturbed elliptic equations.
\newblock {\em Comm. Pure Appl. Math.}, 51(11-12):1445--1490, 1998.

\bibitem{Luo2019}
Peng Luo, Shuangjie Peng, Chunhua Wang, and Chang-Lin Xiang.
\newblock Multi-peak positive solutions to a class of {K}irchhoff equations.
\newblock {\em Proc. Roy. Soc. Edinburgh Sect. A}, 149(4):1097--1122, 2019.

\bibitem{Noussair2000}
Ezzat~S. Noussair and Shusen Yan.
\newblock On positive multipeak solutions of a nonlinear elliptic problem.
\newblock {\em J. London Math. Soc. (2)}, 62(1):213--227, 2000.

\bibitem{Oh1988}
Yong-Geun Oh.
\newblock Existence of semiclassical bound states of nonlinear
  {S}chr\"{o}dinger equations with potentials of the class {$(V)_a$}.
\newblock {\em Comm. Partial Differential Equations}, 13(12):1499--1519, 1988.

\bibitem{Oh1990}
Yong-Geun Oh.
\newblock On positive multi-lump bound states of nonlinear {S}chr\"{o}dinger
  equations under multiple well potential.
\newblock {\em Comm. Math. Phys.}, 131(2):223--253, 1990.

\bibitem{Rabinowitz.1992}
Paul~H Rabinowitz.
\newblock On a class of nonlinear schr{\"o}dinger equations.
\newblock {\em Zeitschrift f{\"u}r Angewandte Mathematik und Physik (ZAMP)},
  43(2):270--291, 1992.

\bibitem{Rey1990}
Olivier Rey.
\newblock The role of the {G}reen's function in a nonlinear elliptic equation
  involving the critical {S}obolev exponent.
\newblock {\em J. Funct. Anal.}, 89(1):1--52, 1990.

\bibitem{Wang2012}
Jun Wang, Lixin Tian, Junxiang Xu, and Fubao Zhang.
\newblock Multiplicity and concentration of positive solutions for a
  {K}irchhoff type problem with critical growth.
\newblock {\em J. Differential Equations}, 253(7):2314--2351, 2012.

\bibitem{Wang1993}
Xuefeng Wang.
\newblock On concentration of positive bound states of nonlinear
  {S}chr\"{o}dinger equations.
\newblock {\em Comm. Math. Phys.}, 153(2):229--244, 1993.

\bibitem{Wang2020}
Zhuangzhuang Wang, Xiaoyu Zeng, and Yimin Zhang.
\newblock Multi-peak solutions of {K}irchhoff equations involving subcritical
  or critical {S}obolev exponents.
\newblock {\em Math. Methods Appl. Sci.}, 43(8):5151--5161, 2020.

\bibitem{Zhang2017}
Jianjun Zhang.
\newblock Standing waves with a critical frequency for nonlinear
  {S}chr\"{o}dinger equations involving critical growth.
\newblock {\em Appl. Math. Lett.}, 63:53--58, 2017.

\bibitem{Zhang2014}
Jianjun Zhang, Zhijie Chen, and Wenming Zou.
\newblock Standing waves for nonlinear {S}chr\"{o}dinger equations involving
  critical growth.
\newblock {\em J. Lond. Math. Soc. (2)}, 90(3):827--844, 2014.

\bibitem{Zhang2015}
Jianjun Zhang and Wenming Zou.
\newblock Solutions concentrating around the saddle points of the potential for
  critical {S}chr\"{o}dinger equations.
\newblock {\em Calc. Var. Partial Differential Equations}, 54(4):4119--4142,
  2015.

\end{thebibliography}

\end{document}